\documentclass[12pt]{article}
\usepackage[american]{babel}
\usepackage{amsmath}
\usepackage{latexsym}
\usepackage{amssymb}
\usepackage{theorem}
\usepackage{diagrams}
\usepackage{mathrsfs}
\usepackage{calligra}
\DeclareMathAlphabet{\mathcalligra}{T1}{calligra}{m}{n}
\DeclareFontShape{T1}{calligra}{m}{n}{<->s*[1.2]callig15}{}
\theoremstyle{change}
\newtheorem{Thm}{Theorem}
\newtheorem{Cor}[Thm]{Corollary}

\newtheorem{Lem}[Thm]{Lemma}
{\theorembodyfont{\rmfamily}
\newtheorem{Num}[Thm]{}

\newtheorem{Def}[Thm]{Definition}}

\renewcommand{\phi}{\varphi}
\renewcommand{\rho}{\varrho}
\newcommand{\Cen}{\mathrm{Cen}}
\newcommand{\fa}{\mathfrak  a}
\newcommand{\fs}{\mathfrak  s}
\newcommand{\fm}{\mathfrak  m}
\newcommand{\fl}{\mathfrak  l}
\newcommand{\fh}{\mathfrak  h}
\newcommand{\fn}{\mathfrak  n}
\newcommand{\fz}{\mathfrak  z}
\newcommand{\Ad}{\mathrm{Ad}}
\newcommand{\ot}{\leftarrow}

\newcommand{\proof}{\par\medskip\rm\emph{Proof. }}

\newcommand{\qed}{\ \hglue 0pt plus 1filll $\Box$}

\newcommand{\mapstoo}{\longmapsto}

\newcommand{\RR}{\mathbb{R}}
\newcommand{\ZZ}{\mathbb{Z}}

\renewcommand{\SS}{\mathbb{S}}

\newcommand{\QQ}{\mathbb{Q}}
\newcommand{\CC}{\mathbb{C}}

\newcommand{\SKIP}[1]{}
\newcommand{\SO}{\mathrm{SO}}

\newcommand{\SU}{\mathrm{SU}}

\newcommand{\AutLie}{\mathrm{AutLie}}
\newcommand{\Lie}{\mathrm{Lie}}
\newcommand{\Gal}{\mathrm{Gal}}

\newcommand{\fk}{\mathfrak{k}}
\newcommand{\fsl}{\mathfrak{sl}}

\newcommand{\Aut}{\mathrm{Aut}}

\newcommand{\Nor}{\mathrm{Nor}}

\newcommand{\inv}{{-1}}

\renewcommand{\setminus}{-}
\begin{document}

\title{{\bf The topology of a semisimple Lie group is essentially unique}}
\author{Linus Kramer\thanks{Supported by SFB 878.}}
\maketitle
\begin{abstract}
We study locally compact group topologies
on simple and semisimple Lie groups. We show that the Lie group topology on such
a group $S$ is very rigid: every 'abstract' isomorphism between
$S$ and a locally compact and $\sigma$-compact group $\Gamma$ is
automatically a homeomorphism, provided that $S$ is absolutely simple.
If $S$ is complex, then noncontinuous field automorphisms of the
complex numbers have to be considered, but that is all.
We obtain similar results for semisimple groups.
\end{abstract}
Abstract isomorphisms between Lie groups have been studied by several
authors. In particular, \'E.~Cartan \cite{Car} and B.~van der Waerden \cite{vdW}
proved that an abstract isomorphism between compact semisimple Lie groups
is automatically continuous. This was generalized by H.~Freudenthal \cite{Freu}
to isomorphisms between absolutely simple real Lie groups.
A model-theoretic proof of Freudenthal's result was given later
by Y.~Peterzil, A.~Pillay and S.~Starchenko \cite{PPS}.
More generally, abstract isomorphisms between simple algebraic groups
were studied by A.~Borel and J.~Tits~\cite{BT}. 

However, all of these results deal with rigidity within the class of
Lie groups. In the present paper, we study a more general problem:
what can be said about abstract isomorphisms between locally compact
groups and simple Lie groups. To put it differently, we want to determine
to what extent the group topology of a simple Lie group is unique.
Some restrictions on the topology are obviously necessary. We have
to exclude the discrete topology, and also all topologies that
are not locally compact (the field of real numbers admits many 
field topologies which are not locally compact). 
A simplified version of our main result can be stated as follows.
For compact $S$, this was proved
by R.~Kallman \cite{Ka}.

\medskip\noindent\textbf{Theorem}
{\em
Let $S$ be an absolutely simple real Lie group. Then the Lie group topology
is the unique locally compact and $\sigma$-compact group topology on $S$.}

\medskip\noindent
More general results are given in Theorems~\ref{CompactThm}, \ref{MainThmReal},
\ref{MainTheoremInComplexCase}, and \ref{SemisimpleMainTheorem}.

\medskip\noindent
I thank Theo Grundh\"ofer, Karl Heinrich Hofmann, Bernhard M\"uhlherr,
Karl-Hermann Neeb, Reiner Salzmann and the referee for helpful
remarks.
\begin{Num}\textbf{Notation}
\label{Generalities}
By a \emph{simple Lie group} we mean a connected centerless
real Lie group $S$ whose Lie algebra $\Lie(S)=\fs$ is simple.
Such a Lie group is simple as an abstract group,
see Salzmann \emph{et al.}  \cite[94.21]{Szm1}, 
and conversely, every
nondiscrete Lie group which is simple as an abstract group is a simple
Lie group in the sense above. We say that $S$ (or $\fs$)
is \emph{absolutely simple} if $\fs\otimes_\RR\CC$ is a complex simple Lie algebra.

By $\AutLie(S)$ we denote the Lie group of all Lie group automorphisms of
$S$. This is the same group as $\Aut_\RR(\fs)$, the group of all
$\RR$-linear Lie algebra automorphisms. The group $S$ acts
faithfully by conjugation on itself (from the left). In this way we
may view it in a natural way as an open subgroup of $\AutLie(S)$.
In fact, $S=\AutLie(S)^\circ$ and $S$ has finite index in $\AutLie(S)$,
see Helgason \cite[IX Thm.~5.4]{Helgason} for the compact (and complex) case
and Murakami \cite[Cor.~2]{Murakami} for the noncompact real case.
In particular, $\AutLie(S)$ is second countable.
\end{Num}
\begin{Def}
We call a Hausdorff space \emph{$\sigma$-compact} if it is a countable
union of compact subsets (some authors require that $\sigma$-compact
spaces are also locally compact, but the present definition is
more suitable for our purposes). Finite products of
$\sigma$-compact spaces are $\sigma$-compact. If
$f:X\rTo Y$ is a continuous map between Hausdorff 
spaces and if $X$ is $\sigma$-compact, then $f(X)$ is
also $\sigma$-compact. Closed subspaces of $\sigma$-compact
spaces are again $\sigma$-compact.
\end{Def}
All topological groups are assumed to be Hausdorff.
We recall the following results about locally compact groups.
\begin{Thm}[Open Mapping Theorem]\label{OpenMap}
Let $\psi:G\rTo H$ be a surjective continuous homomorphism between locally compact
groups. If $H$ is $\sigma$-compact, then $\psi$ is an open map.
\proof
See Hewitt-Ross \cite[II.5.29]{HR} or Stroppel \cite[6.19]{Stroppel}.
\qed
\end{Thm}
\begin{Thm}[Automatic continuity]\label{BorelMap}
Suppose that $G$ is a locally compact group and that $H$ is a $\sigma$-compact
group. Assume that $\psi:G\rTo H$ is is a group homomorphism which is also
a Borel map, i.e. that the
preimage of every open set $U\subseteq H$ is a Borel set.
Then $\psi$ is continuous.

\proof
This is a special case of Hewitt-Ross \cite[V.22.18]{HR}; see also 
Kleppner \cite[Thm.~1]{Kle}.
\qed
\end{Thm}
\begin{Lem}
\label{Lemma1}
Let $S$ be a simple $m$-dimensional Lie group and $G\subseteq \AutLie(S)$ an open subgroup.
Suppose that $C\subseteq G$ is a compact subset that
contains a nonconstant smooth curve. Then there exist
elements $g_0,\ldots,g_m\in G$
such that $g_0Cg_1Cc_2\cdots g_{r-1}Cg_m$ is a compact neighborhood
of the identity.

\proof
Let $c:(-1,1)\rTo C$ be a smooth curve with tangent vector $\dot c(0)=X\neq 0$.
Translating by $c(0)^{-1}$, we may assume that $c(0)=e$ and that
$X\in T_eG=\Lie(G)$. Since $G$ acts via $\Ad$ irreducibly on $\Lie(G)$,
we find elements $h_1,\ldots,h_m\in G$ such that the vectors 
$\Ad(h_1)(X),\ldots,\Ad(h_m)(X)$ span $\Lie(G)$. 
From the inverse function theorem we see that 
$D=h_1 C h_1^\inv h_2 C h_2^\inv\cdots h_m C h_m^\inv$ is a neighborhood of
the identity. 
\qed
\end{Lem}
The next lemma is essentially due to van der Waerden; see \cite[p.~783]{vdW}
and Freudenthal \cite[Satz~8]{Freu}. We denote the commutator by $[a,b]=aba^{-1}b^{-1}$.
\begin{Lem}
\label{Lemma2}
Let $S$ be a simple $m$-dimensional Lie group and $G\subseteq\AutLie(S)$ an open subgroup.
Assume that $D\subseteq G$ is a compact neighborhood of the identity.
Let $U\subseteq G$ be an arbitrary neighborhood of the identity.
Then there exist elements $a_1,\ldots,a_m\in G$
such that $[a_1,D][a_2,D]\cdots [a_m,D]$ is a neighborhood
of the identity which is contained in $U$.

\proof
Let $a\in G\setminus\{ e\}$. Then $\mathrm{Ad}(a)\neq\mathrm{id}$,
so there exists $X\in \Lie(G)$ with $\mathrm{Ad}(a)(X)-X\neq 0$.
Since $S$ acts irreducibly on $\Lie(G)$, we can find elements
$a_1,\ldots,a_m$ in any neighborhood of the identity, and
vectors $X_1,\ldots,X_m\in \Lie(G)$ such that
$\mathrm{Ad}(a_1)(X_1)-X_1,\ldots,\mathrm{Ad}(a_m)(X_m)-X_m$
is a basis of $\Lie(G)$.
It follows readily from the inverse function theorem that 
$[a_1,D]\cdots[a_m,D]$ is a compact neighborhood of the identity.

Let now $U$ be an open neighborhood of the identity.
Consider the continuous map $h:G^m\times G^m\rTo G$,
$(x_1,\ldots,x_m,y_1,\ldots,y_m)\mapstoo
[x_1,y_1][x_2,y_2]\cdots[x_m,y_m]$.
We have that $h(\{e\}^m\times D^m)\subseteq U$. By
Wallace's Lemma, see Kelley \cite[5.12]{Kel}, there is an open neighborhood
$V$ of the identity such that $h(V^m\times D^m)\subseteq U$.
The claim follows if we choose $a_1,\ldots,a_m\in V$.
%
%
\qed
\end{Lem}
The following technical result is the main ingredient in our
continuity proofs. It generalizes Kallman's method \cite{Ka}.
\begin{Thm}
\label{MainProp}
Let $\Gamma$ be a locally compact and $\sigma$-compact group. Let
$S$ be a simple Lie group, let $G$ be an open subgroup of
$\AutLie(S)$ and suppose that 
\[
\phi:\Gamma\rTo G
\]
is an
abstract surjective group homomorphism. Suppose that there is a compact subset
$C\subseteq G$ which contains a nonconstant
smooth curve and whose $\phi$-preimage
$\phi^{-1}(C)$ is $\sigma$-compact.
Then $\phi$ is continuous and open.

\proof
By Lemma~\ref{Lemma1}, there are elements $g_0,\ldots,g_r\in G$
such that $D=g_0Cg_1Cg_2\cdots g_{r-1}Cg_r$ is a compact
neighborhood of the identity. If we choose $\phi$-preimages $g_i'$ of the $g_i$,
then we have $\phi^{-1}(D)=g_0'\phi^{-1}(C)g_1'\phi^{-1}(C)g_2'\cdots g_{r-1}'\phi^{-1}(C)g_r'$.
In particular, $\phi^{-1}(D)$ is $\sigma$-compact.
Let $U\subseteq G$ be an open neighborhood of the identity.
By Lemma~\ref{Lemma2} we find elements $a_1,\ldots a_m\in G$ such that
$E_{a_1,\ldots,a_m}=[a_1,D]\cdots [a_m,D]\subseteq U$ is
a neighborhood of the identity. Moreover, the set
$\phi^{-1}(E_{a_1,\ldots,a_m})$ is $\sigma$-compact and in particular a
Borel set in $\Gamma$.
If $W\subseteq G$ is an arbitrary open subset, then we find
a countable collection of elements $b_j,a_{1,j},\ldots,a_{m,j}\in G$ such that
$W=\bigcup_{j=0}^\infty b_jE_{a_{1,j},\ldots,a_{m,j}}$,
because $G$ is second countable.
Each set $\phi^{-1}(b_jE_{a_{1,j},\ldots,e_{m,j}})$
is Borel, so $\phi^{-1}(W)$ is a Borel set in $\Gamma$.
Therefore $\phi$ is a Borel map and by Theorem~\ref{BorelMap} continuous.
By Theorem~\ref{OpenMap}, $\phi$ is open.
\qed
\end{Thm}
If $S$ is a compact simple Lie group, then every open subgroup $G$ of
$\AutLie(S)$ is compact by the remarks in No.~\ref{Generalities}.
Therefore we have the following consequence
of Theorem~\ref{MainProp}. In the case that $\phi$ is bijective, this result was proved
by Kallman in \cite{Ka}, and for compact $\Gamma$ in Hofmann-Morris \cite[Thm.~5.64]{HofMor}.
\begin{Cor}
\label{CompactThm}
Suppose that $S$ is a compact simple Lie group and that
$G\subseteq \AutLie(S)$ is open. Let $\Gamma$ be a locally
compact and $\sigma$-compact group and assume that
$\phi:\Gamma\rTo G$ is an abstract surjective group homomorphism.
Then $\phi$ is continuous and open.
In particular, $\Aut(S)=\AutLie(S)$.
\qed
\end{Cor}
In order to extend this result to noncompact simple Lie groups, we
need some structure theory. We will avoid the classification
of real simple Lie groups. The following case, however, requires
special considerations in our approach.
\begin{Lem}
\label{Rank1Lemma}
Let $\fs$ be an absolutely simple real Lie algebra of real rank
$\mathrm{rk}_\RR(\fs)=1$.
Let $\fs=\mathfrak k\oplus\fa\oplus\fn$ be an Iwasawa
decomposition and put $\fm=Cen_{\mathfrak k}(\fa)$.
If $\fm$ is abelian, then
$\fs=\mathfrak{sl}_2(\RR)$ or $\fs=\mathfrak{su}_{2,1}(\CC)$
(the Lie algebra of the special unitary group of a $3$-dimensional
nondegenerate complex hermitian form of Witt index $1$).

\proof
We note that $\dim(\fa)=\mathrm{rk}_\RR(\fs)=1$.
If $\fm$ is abelian, then $(\fm\oplus\fa)\otimes_\RR\CC$ is a
Cartan subalgebra in $\fs\otimes_\RR\CC$; see Knapp \cite[6.47]{Knapp} or
Helgason \cite[p.~259]{Helgason}. Thus $\fs$ is quasi-split:
all nodes in the Tits diagram are white/encircled (see Tits \cite{TitsAlg}, and also
the ``Satake diagrams''
in \cite[Table VI]{Helgason} or in Warner \cite{Warner}).
Because $\fs$ has real rank $1$, all white/encircled nodes of the underlying Dynkin diagram
are in one orbit of the $\Gal(\CC/\RR)$-action. Since
$\fs$ is absolutely simple, the underlying Dynkin diagram is
connected and hence a tree. Thus the Dynkin diagram is either $\mathsf A_1$
or $\mathsf A_2$. The corresponding quasi-split real Lie algebras are
$\mathfrak{sl}_2(\RR)$ and $\mathfrak{su}_{2,1}(\CC)$ (this follows
either from the tables in Helgason \cite[p.~259]{Helgason} or 
directly from the classification of involutions on $\fsl_3(\CC)$;
see also Tits \cite[p.~55]{TitsAlg}).
\qed
\end{Lem}
We need the following characterization of real absolutely simple
Lie algebras which we could not find in the literature.
The result follows of course also from the classification of the real
simple Lie algebras, see Helgason \cite[pp.~532--534]{Helgason}.
We remark that the complex Kac Moody algebra of type
$\tilde{\mathsf A}_{2k+1}$ has, for $k\geq 1$, a real form where every rank~$1$ Levi 
factor is of type $\fsl_2(\CC)$, so the result is not completely trivial.
\begin{Lem}
\label{ComplexLemma}
Let $\fs$ be an absolutely simple real Lie algebra. Then there exists a
next-to-minimal parabolic subalgebra whose semisimple Levi algebra
is not isomorphic to $\fsl_2(\CC)$.

\proof
Assume that the semisimple Levi algebra of every next-to-minimal parabolic
subalgebra
is isomorphic to $\fsl_2(\CC)$. Then the underlying Dynkin diagram of each
corresponding complexified Levi algebra is $\mathsf A_1\times\mathsf A_1$
and the Galois group $\Gal(\CC/\RR)$ permutes the two nodes.
In the Tits diagram, the two nodes are white/encircled.

Now in general, the Tits diagram of a semisimple Levi group of a parabolic
is obtained by removing from the Tits diagram of $\fs$
a collection of $\Gal(\CC/\RR)$-orbits  consisting of white/encircled nodes. 
Thus, in our situation all nodes are white/encircled (hence $\fs$ is quasi-split)
and all $\Gal(\CC/\RR)$-orbits consist of two white/encircled nodes which do not
form an edge in the Dynkin diagram. It follows that $\Gal(\CC/\RR)$
acts freely on the Dynkin diagram.

On the other hand, the Dynkin diagram of $\fs\otimes_\RR\CC$ is connected 
(because $\fs$ is absolutely simple) and therefore a tree.
This is a contradiction: the Galois group $\Gal(\CC/\RR)\cong\ZZ/2$
cannot act freely on a tree.
\qed
\end{Lem}
The following theorem contains and extends Freudenthal's Continuity Theorem \cite{Freu}.
\begin{Thm}
\label{MainThmReal}
Suppose that $S$ is an absolutely simple Lie group and that
$G\subseteq \AutLie(S)$ is open. Let $\Gamma$ be a locally compact
and $\sigma$-compact
group and assume that $\phi:\Gamma\rTo G$ is an abstract group isomorphism.
Then $\phi$ is a homeomorphism. In particular, $\Aut(S)=\AutLie(S)$.
\end{Thm}
Before we embark on the proof, we note the following.
We may assume that the real rank $\ell=\mathrm{rk}_\RR(\fs)$ of $S$ is at least $1$, since we
dealt already with compact simple groups (groups of real rank $\ell=0$)
in Corollary~\ref{CompactThm}. We first consider the case $\ell=1$ in a slightly
more general situation.
Suppose that $H$ is a (not necessarily connected)
reductive real Lie group of real rank~$1$
(i.e. the semisimple part of $H^\circ$ has real rank $1$).
Let $\Lie(H)=\fh$ denote its Lie algebra. The semisimple part $\fh_{ss}$ of $\fh$
decomposes as a sum of a simple ideal $\fh_s$ of real rank $1$ and
a compact semisimple ideal $\fh_c$ (which may be trivial).
We assume that there is a
maximal compact subgroup $K\subseteq H$ corresponding to a Cartan involution and 
correspondingly an Iwasawa decomposition
\[
\fh=\fk\oplus\fa\oplus\fn.
\]
From our assumptions, the identity component $H^\circ$ has finite index in $H$.
Let $\fm=\Cen_\fk(\fa)$ and $A=\exp(\fa)$. The compact group $M=\Cen_K(A)$ has $\fm$ as its Lie algebra
and $\fm$ contains $\fh_c$. We distinguish the following cases.

\medskip\noindent
Case (A): \emph{$\fm$ is not abelian}. Then we find an element $h\in M^\circ$
such that the conjugacy class $C=\{ghg^\inv\mid g\in M\}$ is
compact and of positive dimension. Put $L=\Cen_H(A)$. Then $L=MA$ is a central
product and therefore
\[
C=\{ghg^\inv\mid g\in M\}=\{ghg^\inv\mid g\in L\}.
\]

\medskip\noindent
Case (B): \emph{$\fm$ is abelian}.
Then $\fh_c=0$ and thus the semisimple part $\fh_{ss}=\fh_s$ is in fact simple of real rank $1$.
Thus we may use Lemma~\ref{Rank1Lemma}.
If $\fh_s\neq\fsl_2(\CC)$, then $\fk\cap\fh_s$ is isomorphic to
$\RR$ or $\RR\oplus\mathfrak{su}(2)$; in particular, $\fk\cap\fh_s$
has a $1$-dimensional center $\fz$.  Assume that this is the case.
Let $Z\subseteq K$ denote the corresponding connected
subgroup. Since $\fk\cap\fh_{s}$ is an ideal in $\fk$, the group $K^\circ$
centralizes $Z$. Let $h$ be an element in the analytic subgroup $H_s$
corresponding to $\fh_s$ whose $Z$-conjugacy class has positive
dimension. Let $L=\Cen_H(Z)$. The Lie algebra $\Lie(L)=\fl$ decomposes
as $\fl\cong \Cen(\fh)\oplus\fz$. Since $L$ is a finite
extension of $L^\circ$, the set
\[
C=\{ghg^\inv\mid g\in L\}
\]
is compact and of positive dimension.

\medskip\noindent
The remaining case, where $\fh_s=\fsl_2(\CC)$, will not be important.

\medskip
\emph{Proof of Theorem \ref{MainThmReal}.}
We use the structure theory of the (not necessarily connected) group $G$. 
See Warner \cite[p.~85]{Warner} for some some remarks on the nonconnected case.
By Lemma~\ref{ComplexLemma} we can find a next-to-minimal
parabolic $P\subseteq G$ whose semisimple Levi group is not of type $\fsl_2(\CC)$.
Let $H\subseteq P$ denote the reductive Levi group of $P$.
The group $H$ can be written as a centralizer of a torus; see Warner \cite[p.~73]{Warner}.
From this description it is clear that $\phi^\inv(H)$ is closed in $\Gamma$.

We now use our results above about reductive groups of real rank $1$.
In the cases (A) and (B) above, we see from the respective
definitions of the subgroup $L\subseteq H$ that $\phi^\inv(L)$ is also
closed and hence $\sigma$-compact. Therefore $\phi^\inv(C)$ is in both cases $\sigma$-compact.
The claim follows now from Theorem~\ref{MainProp}.
\qed

\medskip\noindent
It may seem that the previous proof with the two cases
(A) and (B) is too complicated. However,
Theorem~\ref{MainThmReal} is false if $G$ happens to be a complex Lie group.
In this case, the general construction of the subset $C$ with
the properties required in Theorem~\ref{MainProp} will, in general, not be possible.
In such a complex Lie group, the subgroup $M^\circ$ that we
used in our construction of $C$ will always be abelian.
Nevertheless, we can prove something in the complex case. Our
methods are, however, somewhat different.
We use the following results.
\begin{Thm}
\label{GrMo}
Let $G$ be a locally compact group. If $G/\Cen(G)$ is compact,
then the algebraic commutator group $DG$ of $G$ has compact closure.

\proof
See Grosser-Moskowitz \cite[Cor.~1, p.~331]{GrMo}.
\qed
\end{Thm}
The following is well-known; actually, we need it only for
the group $G=\SU(2)$, where it is easily verified by hand.
\begin{Thm}
\label{Goto}
Let $G$ be a compact semisimple Lie group. Then $G$ consists of commutators,
$G=\{[a,b]\mid a,b\in G\}$.

\proof
See Hofmann-Morris \cite[Thm.~6.55]{HofMor}.
\qed
\end{Thm}
Finally, we use the following fact about the complex numbers. We recall
that $\CC$ has $2^{2^{\aleph_0}}$ (noncontinuous) field automorphisms.
\begin{Thm}
\label{ComplexNumbersUnique}
Let $\mathcal T$ be a nondiscrete locally compact Hausdorff topology on 
the set $\CC$. Suppose that for every $a\in\CC$, the
the maps $z\mapstoo a+z$ and $z\mapstoo az$ are continuous with
respect to $\mathcal T$.
Then there is a field automorphism $\alpha\in\Aut(\CC)$ such
that $\alpha(\mathcal T)$ is the standard topology on $\CC$.

\proof
By Warner \cite[Thm.~11.17]{SWarner}, $(\CC,\mathcal T)$ is a topological
field. By Weil \cite[I.\S3, Thm.~5]{Weil}, there is, up to topological isomorphism, only one
nondiscrete locally compact algebraically closed field; see also Salzmann \emph{et al.}
\cite[58.8]{Szm2}.
\qed
\end{Thm}
\begin{Num}\textbf{Complex simple Lie groups.}
Suppose that $S$ is a complex simple Lie group with Lie algebra $\Lie(S)=\fs$.
The group $\Aut_\CC(\fs)$ of all $\CC$-linear automorphisms
is a complex Lie group containing $S$. We denote by $\Aut_\QQ(\fs)$
the group of all semilinear automorphisms of $\fs$ (with respect to arbitrary
field automorphisms of $\CC$).

The group $\Aut_\CC(\fs)$ is a complex linear algebraic group. 
It can be realized as a matrix group which is defined by a (finite)
set of polynomial equations on the entries. In this way, one obtains
an action of $\Aut(\CC)$ on $\Aut_\CC(\fs)$ (and on $\fs$), where the
field automorphisms are applied entry-wise to the matrices%
\footnote{This follows also if $S$ is viewed as a group scheme defined over $\QQ$, but the
present down-to-earth approach with matrix groups suffices for our purposes.}.
In particular, there are split short exact sequences
\begin{diagram}[height=2em]
1 & \rTo & \Aut_\CC(\fs) & \rTo & \Aut_\RR(\fs) & \rTo^j_\ot& \Gal(\CC/\RR) & \rTo 1\\
&&\dEq&&\dInto&&\dInto\\
1 & \rTo & \Aut_\CC(\fs) & \rTo & \Aut_\QQ(\fs) & \rTo^j_\ot& \Aut(\CC) & \rTo 1\rlap{.}
\end{diagram}
The group $\Aut_\QQ(\fs)$ is thus contained in the group $\Aut(S)$ of all ``abstract''
automorphisms of $S$.
We will see below that both groups are equal.
For $\alpha\in\Aut_\QQ(\fs)$ we put 
\[
c_\alpha=[g\mapstoo {}^\alpha g=\alpha g\alpha^\inv].
\]
\end{Num}
As in the real case, we first consider groups of rank $1$. We remark that
$\mathrm{PSL}_2(\CC)$ has index $2$ in $\AutLie(\mathrm{PSL}_2(\CC))$;
the quotient is $\Gal(\CC/\RR)$.
\begin{Lem}
\label{ComplexRank1}
Suppose that $S=\mathrm{PSL}_2(\CC)$ and $\fs=\fsl_2(\CC)$. Let 
$G\subseteq \Aut_\QQ(\fs)$ be a subgroup
containing $S$ and let $\Gamma$ be a locally compact and
$\sigma$-compact group. Suppose that $\phi:\Gamma\rTo G$
is an abstract group isomorphism.
Then there exists an element $\alpha\in\Aut_\QQ(\fs)$ such that
$c_\alpha\circ\phi$ is a homeomorphism onto an open subgroup of $\AutLie(S)$. In particular,
$\Aut(S)=\Aut_\QQ(\fs)$.

\proof
We represent the elements of $\mathrm{SL}_2(\CC)$ in the standard way
as complex $2\times 2$ matrices.
This gives us a canonical action of $\Aut(\CC)$ on $S=\mathrm{PSL}_2(\CC)$.
Let $V\subseteq S$ denote the unipotent subgroup represented by all upper
triangular matrices with ones on the diagonal. Then $V=\Cen_G(V)$ 
is isomorphic to the additive group $(\CC,+)$. Let
$T\subseteq S$ denote the group presented by all diagonal matrices.
Then $T=\Cen_G(T)$ acts on $V$ as multiplication by (squares of)
nonzero complex numbers.

The group $V'=\phi^\inv(V)=\Cen_\Gamma(V')$ is a closed and therefore
locally compact and $\sigma$-compact copy of the abstract group $(\CC,+)$.
Moreover, the multiplication by any complex scalar is continuous on $V'$,
since each element of $T'=\phi^{-1}(T)$ acts continuously on $V'$. 
By Theorem~\ref{ComplexNumbersUnique}, there is a
field automorphism $\alpha\in \Aut(\CC)$
such that the restriction $c_\alpha\circ\phi|_{V'}:V'\rTo V$ is a
homeomorphism of topological groups.

We now consider the action of $\Gamma$ on the complex projective line
$\CC\mathrm{P}^1$ via $c_\alpha\circ\phi$.
The $\Gamma$-stabilizer of a suitable point is $\Nor_\Gamma(V')$.
Since $S$ acts transitively on $\CC\mathrm{P}^1$, we have a
factorization $\Gamma=\Nor_\Gamma(V')(c_\alpha\circ\phi)^\inv(S)$.
Thus $\Gamma$ and $\Nor_\Gamma(V')$ have the same images in
$\Aut(\CC)$ under $j\circ c_\alpha\circ\phi$.
But $\Nor_\Gamma(V')$ acts continuously
on $V'\cong\CC$, hence it maps into $\Gal(\CC/\RR)$.
Thus $^\alpha G=c_\alpha(\phi(\Gamma))\subseteq\AutLie(S)$.
Now we can apply Theorem~\ref{MainProp}.
Let $C\subseteq V$ denote the closed unit disk.
Then $(c_\alpha\circ\phi)^\inv(C)$ is closed in $V'$ and hence
closed in $\Gamma$. By Theorem~\ref{MainProp}, 
the composite $c_\alpha\circ\phi:\Gamma{\rTo} {{}^\alpha G}\subseteq\AutLie(S)$ is a
homeomorphism.
\qed
\end{Lem}
\begin{Num}
For groups of higher rank,
we recall a few combinatorial facts. Associated to a complex simple
Lie group $S$ of rank $\ell$ there is an $\ell-1$-dimensional
simplicial complex, the spherical building
$\Delta(S)$. The groups $S$ and $\Aut_\QQ(\fs)$ act on $\Delta$
(the latter by not necessarily type-preserving automorphisms).
The group $\Aut_\QQ(\fs)$ acts on the Dynkin diagram
of $S$. In this action, $\Aut_\QQ(\fs)$ cannot exchange long
and short roots, that is, the action is trivial for the diagrams
$\mathsf C_n$, $\mathsf F_4$ and $\mathsf G_2$. To see this,
we note that $\Aut_\QQ(\fs)$ can be split as a semidirect product
of $\Aut_\CC(\fs)$ and $\Aut(\CC)$, with $\Aut(\CC)$ acting trivially
on the Dynkin diagram. This reduces the claim to the group $\Aut_\CC(\fs)$,
where it is a well-known fact; see Jacobson \cite[IX, Thm.~4]{Jac}.

Pairs of opposite roots determine
walls in $\Delta(S)$; these walls are combinatorial $\ell-2$-spheres.
The result that we will need below is that $\Aut_\QQ(\fs)$
in its action on $\Delta(S)$ has the same orbits on the walls
as $S$: one orbit if the Dynkin diagram is simply laced, and two
orbits for the Dynkin diagrams $\mathsf C_n$, $n\geq 2$,
$\mathsf F_4$ and $\mathsf G_2$.
\end{Num}
\begin{Thm}
\label{MainTheoremInComplexCase}
Let $S$ be a complex simple Lie group with Lie algebra $\fs$,
let $G\subseteq \Aut_\QQ(\fs)$ be
a subgroup containing $S$ and let $\Gamma$ be a locally compact and
$\sigma$-compact group. Suppose that $\phi:\Gamma\rTo G$ is
an abstract group isomorphism. Then there exists an element
$\alpha\in Aut_\QQ(\fs)$ such that the composite
$c_\alpha\circ\phi:\Gamma\rTo \Aut_\QQ(\fs)$ is a homeomorphism onto an open
subgroup of $\AutLie(S)$. In particular, $\Aut(S)=\Aut_\QQ(\fs)$.

\proof
Let $H\subseteq S$ be a reductive Levi subgroup of rank $1$ in a next-to-minimal
parabolic $P\subseteq S$. Thus $H\cong H_0 T$, where
$H_0$ is isomorphic to $\mathrm{PSL}_2(\CC)$ or $\mathrm{SL}_2(\CC)$ and $T\cong (\CC^*)^{\ell -1}$
where $\ell$ is the complex rank of $S$. In view of Lemma~\ref{ComplexRank1}
we may assume that $\ell\geq 2$, so $T$ is nontrivial.
We can arrange the matrix representation of $\Aut_\CC(\fs)$ in such a way that
$T$ is a group of diagonal matrices which is invariant under $\Aut(\CC)$.

Let $L=\Cen_G(T)$. We claim that $L$ is contained in $\Aut_\CC(\fs)$.
The group $\Aut(\CC)$ normalizes $T$, hence
every element $g\in L$ is a product $g=h\eta$, with $\eta\in\Aut(\CC)$
and $h\in\Nor_{\Aut_\CC(\fs)}(T)$. Since $T$ is nontrivial, we
find a nontrivial algebraic character $\lambda:T\rTo\CC^*$ which commutes with
$\Aut(\CC)$ (by evaluating a suitable matrix entry on the diagonal).
Then $h^{-1}$ has to act in the same way 
on $\CC^*$ as $\eta$. However, the only nontrivial algebraic automorphism
of $\CC^*$ is inversion, and this map is not induced by a field automorphism
of $\CC$ (because it is not additive).
It follows that $\eta=1$ and thus $g=h\in \Aut_\CC(\fs)$.

The group $L$ is thus an algebraic finite extension of $H$ and acts 
algebraically on $H/\Cen(H)\cong\mathrm{PSL}_2(\CC)$, with kernel $\Cen(L)$. Thus 
\[
L/\Cen(L)\cong H/\Cen(H)\cong\mathrm{PSL}_2(\CC)
\]
(here we use that $\Aut_\CC(\fsl_2(\CC))=\mathrm{PSL}_2(\CC)$).

Now we consider the preimage
$L'=\phi^{-1}(L)=\Cen_\Gamma(\phi^\inv(T))$. This is a closed subgroup of $\Gamma$.
The map $\phi$ induces an abstract group isomorphism
\[
 \tilde\phi:L'/\Cen(L')\rTo L/\Cen(L)\cong\mathrm{PSL}_2(\CC)
\]
between the locally compact and $\sigma$-compact group $L'/\Cen(L')$ and
the Lie group $L/\Cen(L)$.
By Lemma~\ref{ComplexRank1} there is an element $\alpha\in\Aut(\CC)$
such that $c_\alpha\circ\tilde\phi$ is a homeomorphism.

We claim that $^\alpha G\subseteq\AutLie(S)$. We note that the Levi group
$H\subseteq S$ is the pointwise stabilizer of a unique wall $M\cong\SS^{\ell-2}$
in the spherical building $\Delta(S)$. (The wall is the boundary of the flat subspace
corresponding to the vector part of $\Lie(T)$ in the symmetric space
of $S$.) As we remarked above, $\Aut_\QQ(\fs)$ acts on $\Delta(S)$ (by simplicial,
but not necessarily type-preserving maps). The group $\Aut(\CC)$ fixes this
wall $M$ pointwise and acts by type-preserving automorphisms on $\Delta$.
The group $\Nor_{{}^\alpha G}(T)=\Nor_{{}^\alpha G}(L)$ is precisely
the setwise $^\alpha G$-stabilizer of this wall $M$. 
Now $^\alpha G$ has the same orbits
on the walls of $\Delta$ as $S$, whence
$^\alpha G=\Nor_{{}^\alpha G}(T)S$.
The group $\Nor_\Gamma(L')$ acts continuously on $L'/\Cen(L')$. Pushing this action
forward with $c_\alpha\circ \phi$, we see that $\Nor_{{}^\alpha G}(L)$
acts by homeomorphisms on $L/\Cen(L)$. Thus $\Nor_{{}^\alpha G}(L)$
maps into $\Gal(\CC/\RR)\subseteq\Aut(\CC)$, because $\Aut(\CC)$ acts
faithfully on $L/\Cen(L)$.
It follows from the factorization of $^\alpha G$ above
that $^\alpha G$ also maps into $\Gal(\CC/\RR)$, that is,
$^\alpha G\subseteq\AutLie(S)$.

Now we want to apply Theorem~\ref{MainProp}, and we have to find a set $C\subseteq {}^\alpha G$.
Let $K\subseteq H/\Cen(H)$ be a maximal compact subgroup, $K\cong\SO(3)$.
Let $J'\subseteq L'$
denote the preimage of $K$ under the continuous map
\[
L'\rTo L'/\Cen(L')\rTo^{c_\alpha\circ\tilde\phi} L/\Cen(L)= H/\Cen(H).
\]
We note that as an abstract group, $J'$ is a central product of
a group $C'$ which is isomorphic to $\SO(3)$ or $\SU(2)$ and an abelian group.
Moreover $J'/\Cen(J')$ is (via $c_\alpha\circ\tilde\phi$) homeomorphic
to $\SO(3)$ and thus compact.
By Theorem~\ref{GrMo}, the algebraic commutator group $DJ'$ has compact
closure $\overline{DJ'}$. Algebraically, $C'$ is the set of
commutators $C'=\{[g,h]\mid g,h\in J'\}$ by Theorem~\ref{Goto}.
Thus $C'=\{[g,h]\mid g,h\in\overline{DJ'}\}$ is a compact subset of $\Gamma$.
Now $c_\alpha\circ\phi$ maps $C'$ onto a closed subgroup $C\subseteq {}^\alpha G$.
We now may apply Theorem \ref{MainProp} to the map 
$\Gamma{\rTo^{c_\alpha\circ\phi}}{{}^\alpha G}\subseteq\AutLie(S)$
and conclude that $c_\alpha\circ \phi$ is continuous.
\qed
\end{Thm}
We finally consider semisimple groups. If $S$ is an absolutely simple real Lie group
with Lie algebra~$\fs$, then $\Aut_\QQ(\fs)=\Aut_\RR(\fs)=\AutLie(S)$.
We need the following fact. 
\begin{Num}
Let $S$ be a connected centerless semisimple Lie group with Lie algebra $\fs$.
Let $\fs=\fs_1\oplus\cdots\oplus\fs_n$ be its decomposition into simple ideals,
and $S=S_1\times\cdots\times S_n$ the corresponding factorization into simple groups.
Then $\Aut_\RR(\fs)=\AutLie(S)$ is a semidirect product of a subgroup
$\Pi\subseteq\mathrm{Sym}(n)$ and the direct product
$\Aut_\RR(\fs_1)\times\cdots\times\Aut_\RR(\fs_n)$. The group
$\Pi$ consists of all permutations $\pi$ of $\{1,\ldots,n\}$ that preserve
isomorphy of the simple ideals, i.e. $\fs_i\cong\fs_{\pi(i)}$ for all $i$.
Similarly, the group $\Aut(S)$ of all abstract automorphisms of $S$
decomposes as a semidirect product of the same group $\Pi$ and the direct product
$\Aut(S_1)\times\cdots\times\Aut(S_n)$. 
In particular, there is a split exact sequence
\[
1\rTo\Aut(S_1)\times\cdots\times\Aut(S_n)\rTo\Aut(S)\rTo\Pi\rTo 1.
\]
\end{Num}
\begin{Thm}
\label{SemisimpleMainTheorem}
Let $S$ be a connected centerless semisimple Lie group with Lie algebra
$\Lie(S)=\fs$. Let $S=S_1\times\cdots\times S_n$ denote its 
decomposition into simple factors. Let $G\subseteq\Aut(S)$ be
a subgroup containing $S$. Suppose that $\Gamma$ is a locally compact
and $\sigma$-compact group and that $\phi:\Gamma\rTo G$ is an abstract group
isomorphism. Then there exist elements $\alpha_i\in\Aut(S_i)$
such that $\alpha=\alpha_1\times\cdots\times\alpha_n$ conjugates
$G$ to an open subgroup $^\alpha G\subseteq\AutLie(S)$, and
$c_\alpha\circ\phi:\Gamma{\rTo}{}^\alpha G$ is a homeomorphism.
If $S_i$ is absolutely simple, then $\alpha_i$ may be chosen to be
the identity.

\proof
Let $H_i=\prod_{k\neq i}S_k$. The $\Aut(S)$-centralizer of $H_i$ is
$\Cen_{\Aut(S)}(H_i)=\Aut(S_i)$. Thus we have
$\prod_i\Aut(S_i)=\bigcap_i\Nor_{\Aut(S)}(\Cen_{\Aut(S)}(H_i))$.
From this description it
is clear that the subgroup $\Gamma_0=\phi^\inv\left(\prod_i\Aut(S_i)\right)$
is closed (and open, since it has finite index).
By Theorem~\ref{MainThmReal} and Theorem~\ref{MainTheoremInComplexCase},
we find $\alpha_i\in\Aut(S_i)$ such that 
$\alpha=\alpha_1\times\cdots\times\alpha_n$ conjugates $\phi(\Gamma_0)$
into an open subgroup of $\AutLie(S_1)\times\cdots\times\AutLie(S_n)$,
and such that $c_\alpha\circ\phi:\Gamma_0{\rTo}{}^\alpha G$ is a
homeomorphism onto its image. The identity component $\Gamma^\circ$ of
$\Gamma$ is contained in the open subgroup $\Gamma_0$. Therefore
$c_\alpha\circ\phi$ maps $\Gamma^\circ$ homeomorphically onto $S$. Since 
$\Gamma$ acts by automorphisms on $\Gamma^\circ$, the map
$c_\alpha\circ\phi$ extends continuously to $\Gamma\rTo\AutLie(S)$,
that is, $c_\alpha(\phi(\Gamma))\subseteq\AutLie(S)$.
\qed
\end{Thm}
So far, all the groups that we considered were centerless. It is,
however, easy to see that we have the following consequence of
Theorem~\ref{SemisimpleMainTheorem}.
\begin{Cor}
Let $S$ and $G$ be as in Theorem~\ref{SemisimpleMainTheorem}.
Suppose that $\tilde\Gamma$ is a locally compact
and $\sigma$-compact group and that $\phi:\tilde\Gamma\rTo G$ is a central
surjective homomorphism. If $\Cen(\tilde\Gamma)$ is discrete (for example, finite or
countable), then there exist elements $\alpha_i\in\Aut(S_i)$
such that $\alpha=\alpha_1\times\cdots\times\alpha_n$ conjugates
$G$ to an open subgroup $^\alpha G\subseteq\AutLie(S)$, and
$c_\alpha\circ\phi:\tilde\Gamma{\rTo}{}^\alpha G$ is an open map.
If $S_i$ is absolutely simple, then $\alpha_i$ may be chosen to be
the identity.

\proof
We just note that $\tilde\Gamma\rTo\Gamma=\tilde\Gamma/\Cen(\tilde\Gamma)$
is an open map.
\qed
\end{Cor}

\raggedright
Linus Kramer\\
Mathematisches Institut, 
Universit\"at M\"unster,
Einsteinstr. 62,
48149 M\"unster,
Germany\\
\makeatletter
e-mail: {\tt linus.kramer{@}uni-muenster.de}
\end{document}